\documentclass[a4paper,11pt]{article}

\usepackage{amsmath,amssymb}
\usepackage{amsthm}
\usepackage{palatino}
\usepackage{epsfig}
\usepackage{graphicx}

\usepackage{authblk}

\author[1]{Bela Bauer}
\author[2,1]{Claire Levaillant}

\affil[1]{Station Q, Microsoft Research, Santa Barbara, CA 93106, USA}
\affil[2]{Department of Mathematics, University of California, Santa Barbara, CA 93106, USA}

\title{A new set of generators and a physical interpretation for the $SU(3)$ finite subgroup $D(9,1,1;2,1,1)$}

\begin{document}

\maketitle

\begin{center}
\textbf{Abstract}\end{center}After $100$ years of effort, the classification of all the finite subgroups of $SU(3)$ is yet incomplete. The most recently updated list can be found in \cite{PO2}, where the structure of the series $(C)$ and $(D)$ of $SU(3)$-subgroups is studied. We provide a minimal set of generators for one of these groups which has order $162$. These generators appear up to phase as the image of an irreducible unitary braid group representation issued from the Jones-Kauffman version of $SU(2)$ Chern-Simons theory at level $4$. In light of these new generators, we study the structure of the group in detail and recover the fact that it is isomorphic to the semidirect product $\mathbb{Z}_9\times\mathbb{Z}_3\rtimes S_3$ with respect to conjugation.

\section{Introduction and main result}

\subsection{Definition of the group}

Over the past century, there has been interest in classifying and studying the structure of all the finite subgroups of $SU(3)$ because these subgroups appear in physics. A recent paper by Patrick Otto Ludl \cite{PO2} gives in its introduction a complete chronological review on the matter which we skip in the present paper. As of today, the classification of all the finite subgroups of $SU(3)$ remains incomplete. The original classification attempt dates from $1916$ in a work of Blichfeldt \cite{BL}. Two series of finite subgroups of $SU(3)$ named $(C)$ and $(D)$ are then defined but the complete structure of some of these groups was not studied until $2011$ in \cite{PO2}. Our paper focuses on one particular group of type $D$ of order $162$, namely $D(9,1,1;2,1,1)$. The series $(C)$ and $(D)$ respectively contain subseries $\Delta(3n^2)$ and $\Delta(6n^2)$. For a definition of these series, see for instance \cite{PO1}. It was shown in $2009$ by P.O. Ludl in his thesis \cite{PO1} that any $SU(3)$ subgroup of type $(C)$ can be interpreted as a three-dimensional irreducible representation of $\Delta(3n^2)$. As shown by the same author, the group $D(9,1,1;2,1,1)$ provides a counter-example that the same kind of result does not hold for the $SU(3)$ subgroups of type $(D)$, that is there exists at least one $SU(3)$ subgroup of type $D$ that cannot be interpreted as a three-dimensional irreducible representation of $\Delta(6n^2)$. We recall below the definitions of the series $(C)$ and $(D)$ as in \cite{PO2}.
The series $C(n,a,b)$ for $n$ a positive integer and $a$ and $b$ integers with $0\leq a,b\leq n-1$, is defined as the group generated by the permutation matrix $E$ corresponding to the cycle $(1,3,2)$ of the symmetric group $S_3$ and the diagonal matrix $F(n,a,b)$ with $n$-th roots of unity on the diagonal and determinant $1$.
$$\begin{array}{l}C(n,a,b)\,=\,<E, F(n,a,b)>,\;\text{where}\\\\
E=\begin{pmatrix}
0&1&0\\
0&0&1\\
1&0&0
\end{pmatrix}\;\;\;
\text{and}\;\;\;F(n,a,b)=\begin{pmatrix}e^{\frac{2ia\pi}{n}}&0&0\\0&e^{\frac{2ib\pi}{n}}&0\\
0&0&e^{\frac{2i(-a-b)\pi}{n}}\end{pmatrix}\\\end{array}$$
In order to get the $SU(3)$ subgroup $D(n,a,b;d,r,s)$, with $d$ a positive integer and $r$ and $s$ two integers with $0\leq r,s\leq d-1$, add the extra generator\\
$$\begin{pmatrix} e^{\frac{2ir\pi}{d}}&0&0\\0&0&e^{\frac{2is\pi}{d}}\\0&-e^{\frac{2i(-r-s)\pi}{d}}&0\end{pmatrix}$$
\subsection{The new generators}
The interest of our paper is double. We are able to provide only two generators for $D(9,1,1;2,1,1)$ instead of three. Moreover, our two generators are issued from an irreducible unitary braid group representation $B_4\rightarrow U(3)$ obtained by braiding four anyons of topological charge 2 on a fusion tree of total topological charge $0$ in the Jones-Kauffman version of $SU(2)$ level $4$ Chern-Simons theory \cite{KL}. For some explanations and some terminology associated with this theory, we refer the reader to the excellent exposition in \cite{ZW}. We recover the structure of the group as a semidirect product of the normal abelian group $\mathbb{Z}_9\times\mathbb{Z}_3$ and the symmetric group on three letters $S_3$, where the action is given by conjugation. Our results are summarized below.

\newtheorem{theo}{Theorem}
\begin{theo}
Let $t=\frac{\sqrt{2}}{2}\,e^{\frac{2i\pi}{3}}$ and let $\mathcal{G}$ be the subgroup of $SU(3)$ generated by the matrices $G_1$ and $G_2$, defined as follow.
$$G_1=e^{\frac{i\pi}{9}}\,\begin{pmatrix}  2\bar{t}^2&&\\&2t^2&\\&&-2\bar{t}^2\end{pmatrix}\;\;G_2=e^{\frac{i\pi}{9}}\,\begin{pmatrix} t^2&t&-t^2\\
t&0&t\\-t^2&t&t^2\end{pmatrix}$$ The matrices $G_1$ and $G_2$ have the same eigenvalues and so $G_2$ is unitarily similar to $\begin{pmatrix} e^{\frac{i7\pi}{9}}&&\\&-e^{\frac{i4\pi}{9}}&\\&&e^{\frac{-2i\pi}{9}}\end{pmatrix}$. Both matrices thus have order $18$.
\begin{center}\epsfig{file=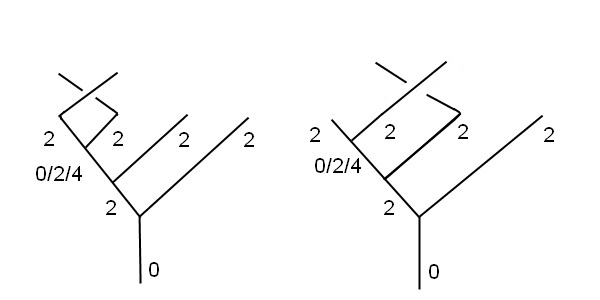, height=5.0cm}\end{center}
\noindent Up to the phase $e^{\frac{i\pi}{9}}$, the matrices $G_1$ and $G_2$ are obtained by respectively braiding the anyons $1$ and $2$ on one hand and $2$ and $3$ on the other hand, like on the above two trees, where these are subject to adequate unitary normalizations. \\
The matrices $G_1$ and $G_2$ satisfy the braid relation $$G_1G_2G_1=G_2G_1G_2$$ and their respective squares commute $$G_1^2\,G_2^2=G_2^2\,G_1^2$$
$\mathcal{G}=\mathcal{N}\overset{\phi}{\rtimes}\mathcal{H}$, where the normal subgroup $\mathcal{N}$ of $\mathcal{G}$ and the subgroup $\mathcal{H}$ of $\mathcal{G}$ are respectively defined by
$$\mathcal{N}=<G_1\,G_2^2\,\,G_1^{-1},\,G_1\,G_2^{-2}\,G_1>\,\simeq\,\mathbb{Z}_9\times\mathbb{Z}_3$$
\begin{multline*}\mathcal{H}=\lbrace I,\,\begin{pmatrix} -1&&\\&-1&\\&&1\end{pmatrix},\, G_1\,G_2\,G_1,\,\begin{pmatrix} -1&&\\&-1&\\&&1\end{pmatrix}\,G_1\,G_2\,G_1\,\begin{pmatrix} -1&&\\&-1&\\&&1\end{pmatrix},\\\,G_1\,G_2\,G_1\,\begin{pmatrix} -1&&\\&-1&\\&&1\end{pmatrix},\,\begin{pmatrix} -1&&\\&-1&\\&&1\end{pmatrix}\,G_1\,G_2\,G_1\rbrace\,\simeq\,S_3\,,\end{multline*}
and where for any matrix $H$ of $\mathcal{H}$, the map $\phi(H)$ is the automorphism of $\mathcal{N}$ which is the conjugation by $H$. Thus $\mathcal{G}$ is isomorphic to $D(9,1,1;2,1,1)$.
\end{theo}

The non-diagonal matrices defining $\mathcal{H}$ are, in the same order as given above,
$$\begin{pmatrix} -\frac{1}{2}& -\frac{1}{\sqrt{2}}&-\frac{1}{2}\\&&\\
-\frac{1}{\sqrt{2}}&0&\frac{1}{\sqrt{2}}\\&&\\
-\frac{1}{2}&\frac{1}{\sqrt{2}}&-\frac{1}{2}\end{pmatrix},\, \begin{pmatrix} -\frac{1}{2}& -\frac{1}{\sqrt{2}}&\frac{1}{2}\\&&\\
-\frac{1}{\sqrt{2}}&0&-\frac{1}{\sqrt{2}}\\&&\\
\frac{1}{2}&-\frac{1}{\sqrt{2}}&-\frac{1}{2}\end{pmatrix},\,\begin{pmatrix}\frac{1}{2}& \frac{1}{\sqrt{2}}&-\frac{1}{2}\\&&\\
\frac{1}{\sqrt{2}}&0&\frac{1}{\sqrt{2}}\\&&\\
\frac{1}{2}&-\frac{1}{\sqrt{2}}&-\frac{1}{2}\end{pmatrix},\,\begin{pmatrix}
\frac{1}{2}& \frac{1}{\sqrt{2}}&\frac{1}{2}\\&&\\
\frac{1}{\sqrt{2}}&0&-\frac{1}{\sqrt{2}}\\&&\\
-\frac{1}{2}&\frac{1}{\sqrt{2}}&-\frac{1}{2}\end{pmatrix}$$

The paper is organized as follows. We define the Hilbert space of the unitary braid group representation whose image is the group that we are studying here and recall some facts in Temperley-Lieb recoupling theory. We then introduce the two generators of the abelian subgroup $\mathcal{N}$ and show that this subgroup is a normal subgroup of $\mathcal{G}=<G_1,G_2>$. We next prove that $\mathcal{G}$ is the semidirect product announced in Theorem $1$.

\section{The unitary braid group representation}

In what follows, $d=-\mathsf{A}^2-\mathsf{A}^{-2}$ is the loop variable in the Temperley-Lieb algebra and the Kauffman variable $\mathsf{A}$ is set to the value $\mathsf{A}=i\,e^{-\frac{2i\pi}{4r}}$, where $k=r-2$ is the level of the theory. Here we work at level $4$ and so $\mathsf{A}=i\,e^{-\frac{i\pi}{12}}$. The set $L=\lbrace 0,1,2,3,4\rbrace$ is the label set of the different particule types. The particules obey fusion rules such that the labels $a$, $b$ and $c$ of a trivalent vertex must satisfy
$$\left\lbrace\begin{array}{l}
a+b+c\;\text{is even}\\
c\leq a+b,\;b\leq a+c,\;a\leq b+c\\
a+b+c\leq 2k
\end{array}\right.$$
Such a triple $(a,b,c)$ is said to be admissible.
By definition, \begin{center}\epsfig{file=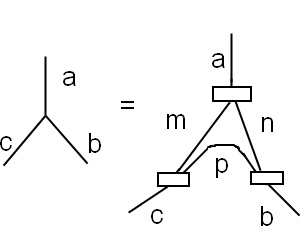, height=5cm}\end{center}
with $$m=\frac{a+c-b}{2},\;n=\frac{a+b-c}{2},\;p=\frac{b+c-a}{2}$$ and where the boxes to the right of the figure represent Jones-Wenzl projectors. We define $\Delta_n$ as the bracket evaluation of the closure of the $(n-1)$-th projector and so $\Delta_1=d$. The classical recursion formula for the Jones-Wenzl projectors shows that $\Delta_n$ is a Chebyschev polynomial of the second kind, whence
$$\Delta_n=(-1)^n\,[n+1]$$ In this formula, $[n]$ denotes the quantum integer taken at $q=\mathsf{A}^2$. $$[n]=\frac{\mathsf{A}^{2n}-\mathsf{A}^{-2n}}{\mathsf{A}^2-\mathsf{A}^{-2}}$$
We have in particular $\Delta_0=\Delta_4=1$ and $\Delta_2=2$.\\
We recall below the $R$-move followed by the $F$-move. The theory appears in \cite{KL}. \begin{center}
\epsfig{file=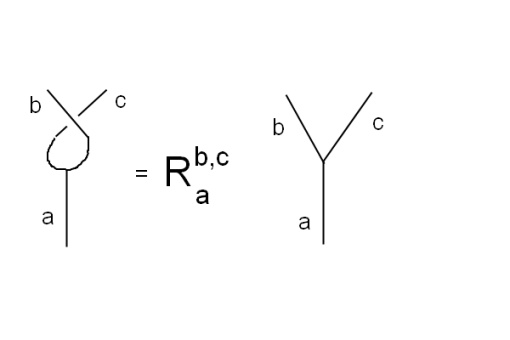, height=6cm}\end{center}
with $$R_a^{b,c}=(-1)^{\frac{b+c-a}{2}}\,\mathsf{A}^{\frac{b(b+c)+c(c+2)-a(a+2)}{2}}$$
In our forthcoming computations, a few useful values are the following. $\overline{R_{0}^{2,2}}=e^{\frac{2i\pi}{3}},\,\overline{R_2^{2,2}}=-e^{\frac{i\pi}{3}},\,\overline{R_4^{2,2}}=e^{-\frac{i\pi}{3}}$, where the bar denotes the complex conjugate.
\begin{center}\epsfig{file=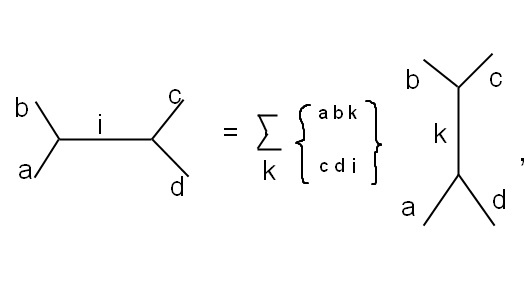, height=5cm}\end{center}
where $k$ runs over all the labels such that both triples $(b,c,k)$ and $(a,d,k)$ are admissible, and where the braces denote the $6j$-symbol \begin{equation}\left\lbrace\begin{array}{ccc}a&b&k\\c&d&i\end{array}\right\rbrace=\frac{T\left[\begin{array}{ccc}
a&b&k\\c&d&i\end{array}\right]\,\Delta_i}{\theta(a,d,i)\theta(c,b,k)}\\
\end{equation}
For the evaluations of the tetrahedral net $T$ and the theta net $\theta$ in terms of the quantum integers, and other formulas from recoupling theory, we refer the reader to the summary in Chapter $9$ of \cite{KL}. We will frequently use that
$$T\left[\begin{array}{ccc} 2&2&j\\2&2&i\end{array}\right]=T\left[\begin{array}{ccc} 2&2&i\\2&2&j\end{array}\right]=\begin{cases} 2&\text{if $(i,j)=(0,0)$},\\
\frac{2}{\sqrt{3}}&\text{if $(i,j)=(2,0)$},\\
0&\text{if $(i,j)=(2,2)$},\\
1&\text{if $(i,j)=(4,0)$},\\
-\frac{1}{\sqrt{3}}&\text{if $(i,j)=(4,2)$},\\
\frac{1}{2}&\text{if $(i,j)=(4,4)$}\end{cases}$$
We consider the $3$-dimensional vector space $V_{4,2,0}$ spanned over $\mathbb{C}$ by the vectors $e_{\alpha}$ with $\alpha\in\lbrace 0,2,4\rbrace$.
\begin{center}\epsfig{file=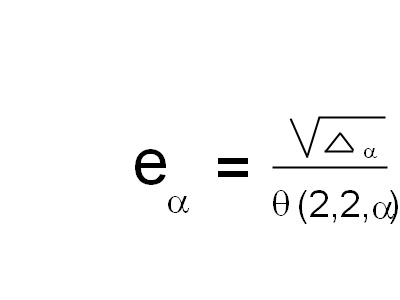, height=3cm}\epsfig{file=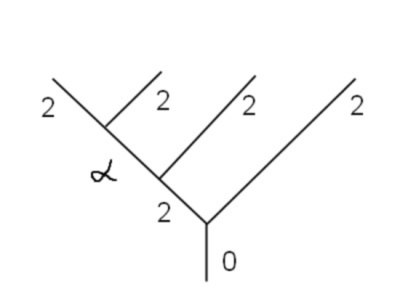, height=4cm}\end{center}
We define an inner product over $V_{4,2,0}$ so that the vectors $e_{\alpha},\,\alpha=0,2,4,$ form an orthonormal basis of $V_{4,2,0}$. The inner product between two such vectors $e_{\beta}$ and $e_{\gamma}$ is obtained by stacking the mirror image of the tree in $e_{\beta}$ on top of the tree in $e_{\gamma}$ and by resolving the crossings with the Kauffman bracket. We now consider the braid group $B_4$ on four strands and we denote by $g_1$, $g_2$ and $g_3$ its three generators. We see that the actions by $g_1$ and $g_3$ on $e_{\alpha}$ can be resolved using an $R$-move and we get in both cases the same matrix
$$\begin{pmatrix}
\overline{R_0^{2,2}}&&\\
&\overline{R_2^{2,2}}&\\
&&\overline{R_4^{2,2}}
\end{pmatrix}$$
 This matrix is the unitary diagonal matrix $G_1$ of Theorem $1$ before the phase $e^{\frac{i\pi}{9}}$ has been added to make it an element of $SU(3)$. Further, the $e_{\beta}$ coordinate of $g_2.e_{\alpha}$ is obtained by evaluating
\begin{center} \epsfig{file=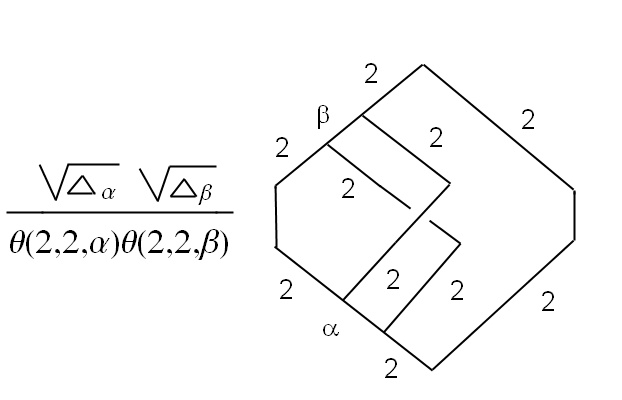, height=7cm}\end{center}
 Since stacking the mirror image of a braid yields the identity, we know that we will obtain a unitary matrix. The braiding in the center of the diagram can be replaced by applying an $F$-move followed by an $R$-move. We obtain
 \begin{center}\epsfig{file=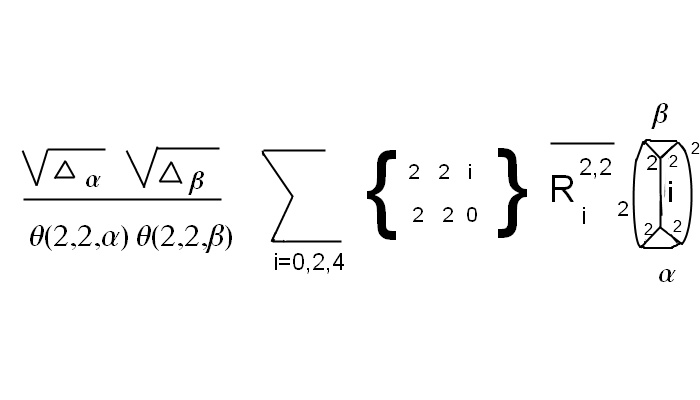, height=5cm}\end{center}
 It remains to compute the network
 \begin{center}\epsfig{file=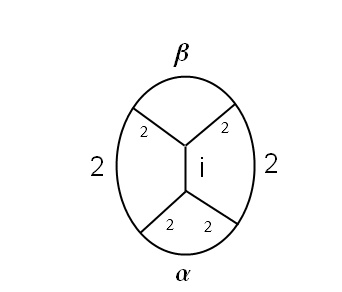, height=5cm}\end{center}
 We do an $F$-move on the lower edge of the figure, labeled $\alpha$. It yields the first member of the equality below. \\
 \epsfig{file=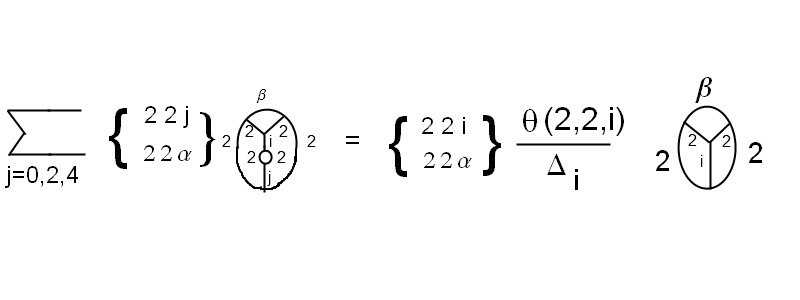, height=5cm}
The figure to the right hand side of this equality is a tetrahedron. After replacing the $6j$-symbol using Eq. $(1)$ and simplifying, we then get
$$\frac{T\left[\begin{array}{ccc}2&2&i\\2&2&\alpha\end{array}\right] T\left[\begin{array}{ccc}2&2&i\\2&2&\beta\end{array}\right]}{\theta(2,2,i)}$$
We get a symmetric matrix whose $(\alpha,\beta)$ coefficient is given by
$$\frac{\sqrt{\Delta_{\alpha}}\sqrt{\Delta_{\beta}}}{\theta(2,2,\alpha)\theta(2,2,\beta)}\,\sum_{i=0,2,4}
\frac{\Delta_i\,\overline{R_i^{2,2}}\;T\left[\begin{array}{ccc}2&2&i\\2&2&\alpha\end{array}\right]\,T\left[
\begin{array}{ccc}2&2&i\\2&2&\beta\end{array}\right]}{\theta(2,2,i)^2}\;,$$
where we used that $T\left[\begin{array}{ccc}2&2&i\\2&2&0\end{array}\right]=\theta(2,2,i)$. Up to the phase $e^{\frac{i\pi}{9}}$ to make it a special unitary matrix, this matrix is the matrix $G_2$ of Theorem $1$. We obtain a unitary braid group representation $B_4\longrightarrow SU(3)$ whose image is $\mathcal{G}=<G_1,G_2>$. \\

In the next part, we show that the subgroup $\mathcal{G}$ of $SU(3)$ generated by the matrices $G_1$ and $G_2$ is isomorphic to the semidirect product $\mathbb{Z}_9\times\mathbb{Z}_3\rtimes S_3$, hence $G_1$ and $G_2$ generate $D(9,1,1;2,1,1)$.
\section{Proof of the Theorem}
\subsection{The abelian subgroup}
We introduce a new orthonormal basis with respect to the inner product defined by stacking a mirror diagram on top of another. \begin{center}
\epsfig{file=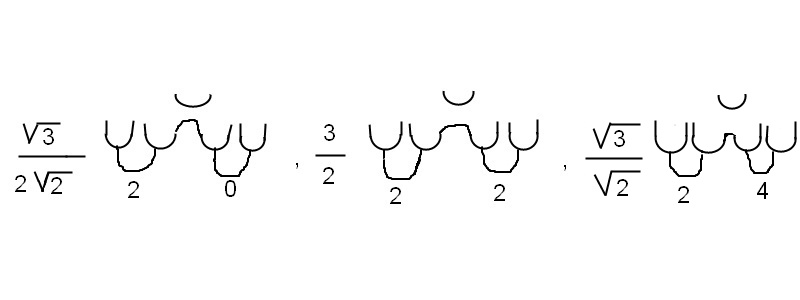, height=6cm}\end{center}
We consider the Hilbert space spanned over $\mathbb{C}$ by these three vectors.
In quantum computing words, the respective normalized diagrams are obtained by taking a quantum trit and an ancilla \begin{center}\epsfig{file=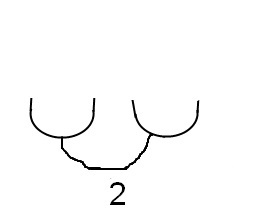, height=3cm}
\epsfig{file=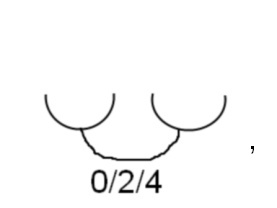, height=3cm}\end{center} by braiding anyons $4$ and $5$ and then doing a forced measurement to zero.
It appears that doing a full twist in the center, that is braiding particules $3$ and $6$ twice leaves the Hilbert space invariant. After doing an additional two braidings between particules $6$ and $7$, and $7$ and $8$, the Hilbert space is still invariant. \\
\begin{center}
\epsfig{file=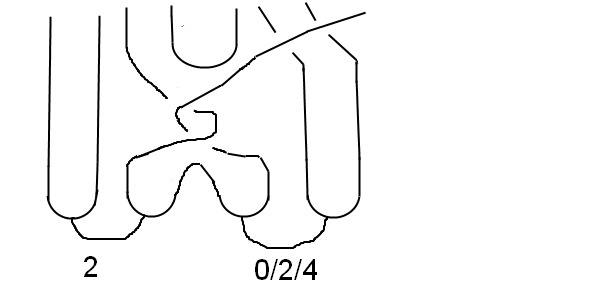, height=6cm}
\end{center}
After computing this action, we obtain the following unitary matrix $$F=\begin{pmatrix}
\frac{-1+i\sqrt{3}}{4}&\frac{\sqrt{2}}{4}(-1+i\sqrt{3})&\frac{1-i\sqrt{3}}{4}\\&&\\\frac{\sqrt{2}}{4}(-1+i\sqrt{3})
&0&\frac{\sqrt{2}}{4}(-1+i\sqrt{3})\\&&\\\frac{-1+i\sqrt{3}}{4}&\frac{\sqrt{2}}{4}(1-i\sqrt{3})&\frac{1-i\sqrt{3}}{4}
\end{pmatrix}$$
This matrix is $G_1G_2\,G_1^{-1}\,G_1^{-1}$. Now consider
$$\begin{array}{l}
A=G_2\,F\,G_2\,F=(G_2F)^2=G_1\,G_2^2\,G_1^{-1},\\
B=G_1\,G_2^{-2}\,G_1
\end{array}$$
and let $\mathcal{N}$ be the subgroup of $\mathcal{G}$ generated by $A$ and $B$.
The matrix $A$ has order $9$ and the matrix $B$ has order $3$. Since
$$B\neq A^3,\,B^2\neq A^3,\,B\neq A^6,\,B^2\neq A^6, $$
we know that $$<A>\cap<B>=\lbrace I\rbrace$$
Moreover, we have $$AB=BA$$
Then, $$\mathcal{N}\simeq\mathbb{Z}_9\times\mathbb{Z}_3$$
We now show that $\mathcal{N}$ is a normal subgroup of $\mathcal{G}$. It suffices to show that for all the integers $k$ and $l$, we have
\begin{eqnarray}
G_1\,A^k\,B^l\,G_1^{-1}\,\in\,<A,B>\\
G_2\,A^k\,B^l\,G_2^{-1}\,\in\,<A,B>
\end{eqnarray}

\newtheorem{Lemma}{Lemma}
\begin{Lemma}
The following are equivalent.
$$\begin{array}{l}(i)\;\forall\,l\,\in\mathbb{Z},\,G_1\,A^l\,G_1^{-1}\,\in\,<A,B>\\ \\ (ii)\;\forall\,(k,l)\,\in\mathbb{Z}^2,\,G_1\,A^l\,B^k\,G_1^{-1}\,\in\,<A,B>\end{array}$$
\end{Lemma}
\noindent\textbf{Proof.} Suppose $(i)$ holds. Write $A^l\,B^k=(AB)^k\,A^{l-k}$. Then,
\begin{eqnarray*}
G_1\,A^l\,B^k\,G_1^{-1}&=& G_1\,G_1^{2k}\,A^{l-k}\,G_1^{-1}\\
&=&G_1^{2k}\,\,G_1\,A^{l-k}\,G_1^{-1}\\
&=&(AB)^k\,\,G_1\,A^{l-k}\,G_1^{-1}\qquad\qquad\qquad\qquad\qquad\qquad\qquad\qquad\qquad\square
\end{eqnarray*}
Thus, in order to prove $(2)$, it suffices to show that $G_1\,A\,G_1^{-1}\in\,<A,B>$. Since by definition,
$A=G_1\,G_2^2\,G_1^{-1}$ and since $G_1^2\,G_2^2=G_2^2\,G_1^2$, we have
$$G_1\,A\,G_1^{-1}=G_2^2$$
As $G_2^2=A^7\,B^2$, point $(2)$ is proven. To have a normal subgroup, it remains to show $(3)$.
Since $G_2^2\,\in\,<A,B>$, we have
$$\begin{array}{cccc}
G_2\,A^k\,B^l\,G_2^{-1}\,\in\,<A,B>&\Longleftrightarrow&G_2^2\,G_2\,A^k\,B^l\,G_2^{-1}\,\in\,<A,B>\\
&&&\\
&\Longleftrightarrow&G_2\,A^{k-2}\,B^{l-1}\,G_2^{-1}\,\in\,<A,B>&(\clubsuit)\\&&&\\
&\Longleftrightarrow&G_2^{-2}\,G_2\,A^k\,B^l\,G_2^{-1}\,\in\,<A,B>\\&&&\\
&\Longleftrightarrow&G_2\,A^{k+2}\,B^{l+1}\,G_2^{-1}\,\in\,<A,B>&(\spadesuit)
\end{array}$$
Equivalences $(\clubsuit)$ and $(\spadesuit)$ show that $G_2\,A^k\,B^l\,G_2^{-1}$ belongs to $<A,B>$ for all $(l,k)\in\mathbb{Z}^2$ is equivalent to $G_2\,A^t\,G_2^{-1}$ belongs to $<A,B>$ for every $t\in\mathbb{Z}$. Thus, it suffices to show that $G_2\,A\,G_2^{-1}\in\,<A,B>$. We have
\begin{eqnarray}
G_2\,A\,G_2^{-1}&=&G_2\,G_1\,G_2\,G_2\,G_1^{-1}\,G_2^{-1}\\
&=&G_1\,G_2\,G_1\,G_2\,G_1^{-1}\,G_2^{-1}\\
&=&G_1^2\,G_2\,G_1\,G_1^{-1}\,G_2^{-1}\\
&=&G_1^2\\
&=&A\,B
\end{eqnarray}
Eq. $(4)$ is by definition of $A$ and Eq. $(5)$ and $(6)$ are obtained by using the braid relation. Hence, we are done with the proof of $(3)$. So, $\mathcal{N}\lhd\mathcal{G}$.

\subsection{The symmetric group}
We now describe the subgroup of $<G_1,G_2>$ that is isomorphic to $S_3$. The product $T_1=G_1G_2G_1$ provides an element of order $2$ in the group, which we take as one of the two generators of the symmetric group $S_3$. Since $G_1$ has order $18$, the product $T_2=G_2\,G_1^9\,G_2^{-1}$ provides another element of order $2$ in the group. And so does $G_2\,G_1^2$. Now set
$$T_3=T_2\,(G_2\,G_1^2)\,T_2=\begin{pmatrix} -1&&\\&-1&\\&&1\end{pmatrix}$$ and set
$$\mathcal{H}=<T_1,T_3>$$
$T_1$ and $T_3$ generate a subgroup $\mathcal{H}$ of $\mathcal{G}$ isomorphic to $S_3$, whose two elements of order $3$ are $T_1\,T_3$ and $T_3\,T_1$. We checked using Mathematica that these two elements of order $3$ are not one of $$A^3,\,A^6,\,A^3\,B,\,A^6\,B,A^3\,B^2,\,A^6\,B^2,\,B,\,B^2$$ Then we must have
$$\mathcal{H}\cap\mathcal{N}=\lbrace I\rbrace$$
\subsection{The semidirect product}
We show that both $G_1$ and $G_2$ can be expressed as the product of an element of $\mathcal{N}$ and an element of $\mathcal{H}$. This will suffice to show that $\mathcal{G}\simeq \mathcal{N}\overset{\phi}{\rtimes}\mathcal{H}$, where
$$\begin{array}{cccc}
\phi:&\mathcal{H}&\longrightarrow& Aut(\mathcal{N})\\
&H&\longmapsto&\left(\begin{array}{ccc}\mathcal{N}&\rightarrow&\mathcal{N}\\
N&\mapsto&HNH^{-1}\end{array}\right)
\end{array}$$
To show that $G_1\in\,\mathcal{N}.\mathcal{H}$, it suffices to show that $G_2^2\,G_1\in\,\mathcal{N}.\mathcal{H}$.
%We have
%\begin{eqnarray}
%G_2^2\,G_1&=&A^{-2}\,B^{-1}\,G_1\\
%&=&G_1\,G_2^{-4}\,G_1^{-2}\,G_2^2\\
%&=&G_1^{-1}\,G_2^{-2}
%\end{eqnarray}
%Eq. $(9)$ is obtained by using that $G_2^2=A^7\,B^2$ as previously used. Eq. $(10)$ is by definition of $A$ %and $B$. Eq. $(11)$ follows from Eq. $(10)$ after using the fact that $G_1^2$ commutes to $G_2^2$.
From $(G_2\,G_1\,G_2)^2=I$ and the braid relation, we derive $(G_2\,G_1)^3=I$. Then,
$$(G_2^2\,G_1)^2=\,G_2\,G_2\,G_1\,G_2\,G_2\,G_1=(G_2\,G_1)^3=I$$
So, $G_2^2\,G_1$ is an element of order $2$ in $\mathcal{G}$. Moreover, $G_2^2\,G_1$ is not one of the three elements of order $2$ of $\mathcal{H}$. Recall that in $(\mathcal{N}\overset{\phi}{\rtimes}\mathcal{H},\star)$, we have
$$(N,H)\star(N,H)=(N\,H\,N\,H^{-1},\,H^2)$$
So, if an element $(N,H)$ has order $2$, then $H$ has order $2$ in $\mathcal{H}$ and $N\,H\,N=H$.
Suppose we can write $G_2^2\,G_1=N\,H$ with $N\,\in\mathcal{N}$ and $H\in\mathcal{H}$. Impose the extra condition $N\,H\,N=H$. Then,
$$G_2^2\,G_1\,N=N\,H\,N=H$$
This implies $N=(G_2^2\,G_1)^{-1}\,H$. Since $\mathcal{N}$ does not contain any element of order $2$, we see that $H$ cannot be the identity matrix. Assume further that $H$ has order $2$. Then, $H$ must be one of
$$T_1\;\text{or}\;T_3\;\text{or}\;T_3\,T_1\,T_3$$
We found out that in the three cases, $(G_2^2\,G_1)^{-1}\,H$ has order $3$. Out of the three possibilities, only $$(G_2^2\,G_1)^{-1}\,T_3$$ belongs to $\mathcal{N}$ and is $A^3\,B$. We get
$$G_2^2\,G_1=A^3\,B\,T_3\;\in\;\mathcal{N}.\mathcal{H}$$
And we conclude
$$G_1=G_2^{-2}\,G_2^2\,G_1\,\in\,\mathcal{N}.\mathcal{H}$$
Similarly, we show that $G_1^2\,G_2\in\mathcal{N}.\mathcal{H}$. Again, since $(G_1\,G_2)^3=I$, this element has order $2$. Also, it does not belong to $\mathcal{H}$. We then use the same strategy as before.
If $G_1^2\,G_2=NH$ with $N\in\mathcal{N},\,H\in\mathcal{H}$ and $N\,H\,N=H$, then
$$N=(G_1^2\,G_2)^{-1}\,H$$
With $H$ of order $2$, this time we find
$$N=(G_1^2\,G_2)^{-1}\,T_3\,T_1\,T_3=B^2$$
And so, $$G_1^2\,G_2=B^2\,T_3\,T_1\,T_3\,\in\mathcal{N}.\mathcal{H}$$
Thus, since $G_1^2=AB$, we get $$G_2=G_1^{-2}\,G_1^2\,G_2\in\mathcal{N}.\mathcal{H}$$
We are ready to conclude. Consider the map
$$\begin{array}{ccc}
\mathcal{N}\overset{\phi}{\rtimes}\mathcal{H}&\overset{\psi}{\longrightarrow}&\mathcal{G}\\
(N,H)&\longmapsto&NH
\end{array}$$
We have seen at the end of paragraph $3.2$ that $\mathcal{H}\cap\mathcal{N}=\lbrace I\rbrace$. Hence $\psi$ is an injective morphism of groups.
By the discussion above, we have 
\begin{eqnarray*}
G_1&=&\psi(A^5\,B^2,\,T_3)\\
G_2&=&\psi(A^{-1}\,B,\, T_3T_1T_3)
\end{eqnarray*}
It follows that
$$<G_1,G_2>\,\subseteq\, Im\,\psi\,\subseteq\, <G_1,G_2>$$
Hence the image of $\psi$ is the whole group $\mathcal{G}$ and
$$\mathcal{G}\simeq\mathcal{N}\overset{\phi}{\rtimes}\mathcal{H}$$

\noindent We end this paper by giving a presentation for the group $\mathcal{G}$. We have
$$
\mathcal{G}=\left<\begin{array}{ccc}A,B,T_1,T_3&|&\begin{array}{l}A^9=B^3=T_1^2=T_3^2=(T_1\,T_3)^3=(T_3\,T_1)^3=I\\
\begin{array}{cc}T_1\,A\,T_1^{-1}=A, &T_3\,A\,T_3^{-1}=A^7\,B^2\\
T_1\,B\,T_1^{-1}=A^6\,B^2,& T_3\,B\,T_3^{-1}=A^3\,B^2\end{array}\end{array}
\end{array}\right>$$
\noindent\textbf{Acknowledgments.} The second author of the paper thanks Zhenghan Wang for suggesting to her the study of the Jones representation from the main theorem and asking her to determine its image. She also thanks him for teaching her some basic recoupling theory and thanks the whole team of Microsoft Research Station Q, where this work was achieved, for great hospitality. We are quite happy to thank Michael Freedman for helpful and enlightening discussions. We thank Patrick Otto Ludl for an informative correspondence.

\end{document}